\documentclass[10pt]{article}
\usepackage{mathrsfs}
\usepackage{amsthm}
\usepackage{amsmath}
\usepackage{graphicx}
\usepackage{color}
\usepackage{amsfonts}
\newtheorem{theorem}{Theorem}[section]

\newtheorem{lemma}[theorem]{Lemma}

\numberwithin{equation}{section}
\normalsize

\begin{document}
\title{\textbf{Fluid Limit of Threshold Voter Models on Tori}}
\author{Xiaofeng Xue \thanks{\textbf{E-mail}: masonxuexf@math.pku.edu.cn \textbf{Address}: School of Mathematical Sciences, Peking University, Beijing 100871, China.}\\ Peking University}

\date{}
\maketitle

\noindent {\bf Abstract}
In this paper, we are concerned with threshold voter models on tori. Assuming that the initial distribution of the process is product measure with density $p$, we obtain a fluid limit of the proportion of vertices in state $1$ as the dimension of the torus grows to infinity. The fluid limit performs a phase transition phenomenon from $p<1/2$ to $p>1/2$.


\section{Introduction}
In this paper, we are concerned with threshold voter models on tori. For integers $d$ and $r$, we denote by $T^d(r)$ the $d$-dimensional torus $\{1,2,\ldots,r\}^d$. In details, for any $x=(x_1,x_2,\ldots,x_i,\ldots,x_d)\in \{1,2,\ldots,r\}^d$ and $1\leq i\leq d$, we define $x^i\in \{1,2,\ldots,r\}^d$ as
\begin{equation}
x^i_j=
\begin{cases}
x_j \text{~if $j\neq i$},\\
x_i+1\text{~if $j=i$ and $x_i<r$},\\
1\text{~if $j=i$ and $x_i=r$}.
\end{cases}
\end{equation}
On $T^d(r)$, for any $x\in \{1,2,\ldots,r\}^d$ and $1\leq i\leq d$, there is an edge connecting $x$ and $x^i$. Therefore, $T^d(r)$ is a regular graph that each vertex has degree $2d$.

Threshold voter model on $T^d(r)$ is with state space $\{0,1\}^{T^d(r)}$, which means that at each vertex there is a spin taking value $0$ or $1$. For any $\eta\in \{0,1\}^{T^d(r)}$ and $x\in T^d(r)$, we denote by $\eta(x)$ the value of $x$. For any $t\geq0$, we denote by $\eta_t$ the configuration of the threshold voter model at moment $t$. For any $t>0$, we define
\[
\eta_{t-}=\lim_{s\uparrow t}\eta_s
\]
as the configuration of the process at the moment just before $t$.
For any $x,y\in T^d(r)$, we say they are neighbors if there is an edge connecting them, denoted by $x\sim y$.

Now we explain how the process evolves. At moment $0$, each vertex of $T^d(r)$ takes $0$ or $1$ according to some probability distribution. Then, the process evolves according to independent Poisson processes $\{N_x(t):t\geq 0\}_{x\in T^d(r)}$. For each $x$, $N_x$ is with rate $1$. The value of $x$ may flip only at event times of $N_x$. For any event time $s$ of $N_x$, if at $s$ there are not less than $d$ neighbors of $x$ taking a different value than $\eta_{s-}(x)$, then at the same moment the value of $x$ flips from $\eta_{s-}(x)$ to $\eta_s(x)=1-\eta_{s-}(x)$, otherwise $\eta_s(x)=\eta_{s-}(x)$. Therefore, the process is a spin system with flip rates given by
\begin{equation}\label{equ of flip rate of threshold voter on tori}
c(x,\eta)=
\begin{cases}
1\text{~if~}\sum_{y:y\sim x}1_{\{\eta(y)\neq\eta(x)\}}\geq d,\\
0\text{~else}.
\end{cases}
\end{equation}
The third section of \cite{LIG1985} gives a precise introduction of spin systems.

Intuitively, $0$ and $1$ are two candidates of an election. Vertices taking $0$ or $1$ are respective supporters of $0$ or $1$. A vertex may change his choice when and only when more than half of the neighbors holding a different opinion. That's why this model is called threshold voter model.

Threshold voter models are introduced by Cox and Durrett in \cite{Dur1991}, where the threshold is assumed to be one. In \cite{Dur1991}, Cox and Durrett gives an important additive dual process of threshold-one voter models and prove that the threshold-one voter model initially with product measure with rate $1/2$ converges weakly to a stationary measure $\nu_{1/2}$. In \cite{Han1999}, Handjani proves a complete convergence theorem for threshold-one voter models on lattices. She shows that with whatever initial distribution, the process converges weakly to a convex combination of three stationary measures $\delta_0$, $\delta_1$ and $\nu_{1/2}$. For threshold voter models on lattices with threshold \(K\geq1\), Liggett and his partners did a lot of important work on judging whether fixation (the process trapped in a state), clustering (all the sites take the same value) or coexistence occurs. The explicit results can be referred in \cite{LIG1992}, \cite{LIG1994} and \cite{LIG1999}. In \cite{Xue2012}, Xue study threshold voter models on homogeneous trees and prove that the critical density of the model is approximately to the proportion of the threshold to the degree of the tree.

Lanchier introduces an opinion dynamics model in \cite{Lanchier2010}, where each vertex is with an opinion in $[0,1]$ and a vertex can mimic a neighbor only if the opinion distance between the vertex pair does not exceed a threshold. The proof of a crucial lemma in this paper is inspired a lot by the approach introduced in \cite{Lanchier2010}.

\section{Main results}
In this section, we give our main results. We consider that the initial distribution of the process is product measure. We aim to describe the proportion of the vertices in state $1$ at any moment $t>0$. Of course the number of vertices in state $1$ is random at any moment, but we will show that as the dimension of the torus grows to infinity the proportion of vertices in state $1$ converges to a deterministic process, which is called the fluid limit.

First we state some notations. For any $A\subseteq T^d(r)$, we denote by $|A|$ the cardinal number of $A$. For any $p\in[0,1]$, we denote by $\mu_p$ the product measure on $\{0,1\}^{T^d(r)}$ with density $p$. In detail, for any $A\subseteq T^d(r)$,
\[
\mu_p\big(\eta:\eta(x)=1,\forall x\in A\big)=p^{|A|}.
\]
In latter sections we assume that $r$ is a fixed integer and not smaller than $2$. For any $t>0$, we denote by $\eta_t^{d,p}$ the configuration at moment $t$ of the threshold voter model on $T^d(r)$ with initial distribution $\mu_p$.  Furthermore, we denote by
\[
A_t^{d,p}=\{x\in T^d(r):\eta_t^{d,p}(x)=1\}
\]
the set of vertices in state $1$ at moment $t$.
The following theorem is our main result about the fluid limit of the process.
\begin{theorem}\label{thm of fluid limit of thres vot on tori}
Assume that $r\geq 2$, then for any $T>0$,
\begin{equation}\label{equ of fluid subcritical}
\lim\limits_{d\rightarrow+\infty}\sup_{0\leq t\leq T}\big|\frac{|A_t^{d,p}|}{r^d}-pe^{-t}\big|=0
\end{equation}
in probability when $p\in [0,1/2)$ and
\begin{equation}\label{equ of fluid supcritical}
\lim\limits_{d\rightarrow+\infty}\sup_{0\leq t\leq T}\big|\frac{|A_t^{d,p}|}{r^d}-1+(1-p)e^{-t}\big|=0
\end{equation}
in probability when $p\in(1/2,1]$.
\end{theorem}

According to the symmetry of $0$ and $1$ in the voter model, it is easy to see that \eqref{equ of fluid supcritical} is a direct corollary of \eqref{equ of fluid subcritical}. In latter sections we will only prove \eqref{equ of fluid subcritical}.

Theorem \ref{thm of fluid limit of thres vot on tori} shows that when $p<1/2$, the proportion of vertices in state $1$ has fluid limit $pe^{-t}$, which converges to $0$ as $t$ grows to infinity. While when $p>1/2$, the fluid limit turns into $1-(1-p)e^{-t}$, which converges to $1$ as $t$ grows to infinity. When $p=1/2$, according to the symmetry of $0$ and $1$, it is easy to see that
\[
P(\eta_t^{d,1/2}(x)=1)=1/2
\]
for any $t>0$ and $x\in T^d(r)$. As a result,
\[
{\rm E}\frac{|A_t^{d,1/2}|}{r^d}\equiv1/2.
\]
Therefore, the threshold voter model performs a phase transition phenomenon from $p<1/2$ to $p>1/2$.

To prove \eqref{equ of fluid subcritical}, we only need to show that for any $\epsilon>0$,
\begin{equation}\label{equ of liminf of subcri fluid}
\lim\limits_{d\rightarrow+\infty}P\big(\inf_{0\leq t\leq T}\{\frac{|A_t^{d,p}|}{r^d}-pe^{-t}\}<-\epsilon\big)=0.
\end{equation}
and
\begin{equation}\label{equ of limsup of subcri fluid}
\lim\limits_{d\rightarrow+\infty}P\big(\sup_{0\leq t\leq T}\{\frac{|A_t^{d,p}|}{r^d}-pe^{-t}\}>\epsilon\big)=0,
\end{equation}
We will prove \eqref{equ of liminf of subcri fluid} and \eqref{equ of limsup of subcri fluid} in Section \ref{sect of proof 2.3} and \ref{sect of proof 2.4}.

\section{Proof of \eqref{equ of liminf of subcri fluid}}\label{sect of proof 2.3}
To prove \eqref{equ of liminf of subcri fluid}, we introduce a Markov process $\{\zeta_t\}_{t\geq 0}$ with state space $\{0,1\}^{T^d(r)}$ to bound $|A_t^{d,p}|$ from below. For any $x\in T^d(r)$,
if $\zeta_0(x)=0$, then $x$ is frozen in state $0$ forever. If $\zeta_0(x)=1$, then $x$ waits for an exponential time $T_x$ with rate $1$ to flip to $0$ and be frozen in $0$ forever. $\{T_x\}_{x\in T^d(r)}$ are independent. In other words, $\{\zeta_t\}_{t\geq 0}$ is a spin system with flip rates given by
\begin{equation}\label{equ of flip rate of zeta}
\widehat{c}(x,\zeta)=
\begin{cases}
1\text{~if $\zeta(x)=1$,} \\
0\text{~if $\zeta(x)=0$.}
\end{cases}
\end{equation}

We write $\zeta_t$ as $\zeta_t^{d,p}$ when $\zeta_0$ is with distribution $\mu_p$ on $T^d(r)$. We denote by
\[
G_t^{d,p}=\{x\in T^d(r):\zeta_t^{d,p}(x)=1\}
\]
the set of vertices in state $1$ at moment $t$.

The following theorem shows the connection between $\zeta_t$ and the threshold voter model $\eta_t$.
\begin{theorem}\label{theorem of couple of eta and zeta}
There exists a couple of $\{\eta_t^{d,p}\}_{t\geq 0}$ and $\{\zeta_t^{d,p}\}_{t\geq 0}$ such that for any $t\geq 0$,
\[
G^{d,p}_t\subseteq A^{d,p}_t.
\]
\end{theorem}
\proof
For any $\eta,\zeta\in T^d(r)$, we write $\eta\gg\zeta$ if and only if $\eta(x)\geq\zeta(x)$ for any $x\in T^d(r)$.
According to \eqref{equ of flip rate of threshold voter on tori} and \eqref{equ of flip rate of zeta}, for any $\eta\gg\zeta, x\in T^d(r)$,
\begin{equation}
\begin{cases}
c(x,\eta)\geq\widehat{c}(x,\zeta)\text{~if $\eta(x)=\zeta(x)=0$,}\\
c(x,\eta)\leq\widehat{c}(x,\zeta)\text{~if $\eta(x)=\zeta(x)=1$.}
\end{cases}
\end{equation}
As a result, Theorem \ref{theorem of couple of eta and zeta} holds according to Theorem 1.5 in Chapter 3 of \cite{LIG1985}.

\qed

Now we construct a martingale about $\zeta_t$ which is important for the proof of \eqref{equ of liminf of subcri fluid}.
\begin{theorem}\label{thm of zeta near a martingle}
For any $t>0$, let $\mathcal{F}_t=\sigma\{\zeta_s^{d,p}:s\leq t\}$. Then, $\{\frac{|G_t^{d,p}|}{pe^{-t}}\}_{t\geq 0}$ is a martingale relative to
$\{\mathcal{F}_t:t\geq 0\}$.
\end{theorem}
\proof
Conditioned on $G_0^{d,p}=A$,
\begin{align}\label{equ (3.3)}
{\rm E}|G_t^{d,p}|&={\rm E}_A|G_t^{d,p}|\notag\\
&={\rm E}\sum_{x\in A}1_{\{\zeta_t^{d,p}(x)=1\}}=\sum_{x\in A}P(\zeta_t^{d,p}(x)=1)\notag\\
&=\sum_{x\in A}P(T_x>t)=|A|e^{-t}
\end{align}
for any $t>0$.

According to Markov property, for any $t>s>0$,
\[
{\rm E}(\frac{|G_t^{d,p}|}{pe^{-t}}\big|\mathcal{F}_s)=\frac{1}{pe^{-t}}{\rm E}_{G_s^{d,p}}|G_{t-s}^{d,p}|.
\]
By \eqref{equ (3.3)},
\[
{\rm E}_{G_s^{d,p}}|G_{t-s}^{d,p}|=|G_s^{d,p}|e^{s-t}.
\]
Therefore,
\[
{\rm E}(\frac{|G_t^{d,p}|}{pe^{-t}}\big|\mathcal{F}_s)=\frac{|G_s^{d,p}|}{pe^{-s}}
\]
and the proof completes.

\qed

Now we give the proof of \eqref{equ of liminf of subcri fluid}.

\proof[Proof of \eqref{equ of liminf of subcri fluid}]
According to Theorem \ref{theorem of couple of eta and zeta},
\begin{align*}
P\big(\inf_{0\leq t\leq T}\{\frac{|A_t^{d,p}|}{r^d}-pe^{-t}\}<-\epsilon\big)
&\leq P\big(\inf_{0\leq t\leq T}\{\frac{|G_t^{d,p}|}{r^d}-pe^{-t}\}<-\epsilon\big)\\
&\leq P\big(\sup_{0\leq t\leq T}[\frac{|G_t^{d,p}|}{pe^{-t}}-r^d]^2>\frac{\epsilon^2}{p^2}r^{2d}\big).
\end{align*}
According to Theorem \ref{thm of zeta near a martingle}, $\frac{|G_t^{d,p}|}{pe^{-t}}$ is a right-continuous martingale and therefore
\[
[\frac{|G_t^{d,p}|}{pe^{-t}}-r^d]^2
\]
is a positive right-continuous submartingale. Then according to Doob inequality,
\begin{align*}
P\big(\sup_{0\leq t\leq T}[\frac{|G_t^{d,p}|}{pe^{-t}}-r^d]^2>\frac{\epsilon^2}{p^2}r^{2d}\big)
&\leq \frac{p^2}{\epsilon^2r^{2d}}{\rm E}[\frac{|G_T^{d,p}|}{pe^{-T}}-r^d]^2\\
&=\frac{e^{2T}}{\epsilon^2r^{2d}}{\rm E}[|G_T^{d,p}|-r^dpe^{-T}]^2.
\end{align*}
Notice that,
\begin{align*}
{\rm E}|G_T^{d,p}|&=\sum_{x\in T^d(r)}P\big(\zeta_T^{d,p}(x)=1\big)\\
&=\sum_{x\in T^d(r)}P\big(\zeta_0^{d,p}(x)=1, T_x>T\big)\\
&=\sum_{x\in T^d(r)}pe^{-T}=r^dpe^{-T}.
\end{align*}
Hence,
\begin{align*}
{\rm E}[|G_T^{d,p}|-r^dpe^{-T}]^2&={\rm Var}(|G_T^{d,p}|)\\
&={\rm Var}(\sum_{x\in T^d(r)}1_{\{\zeta_T^{d,p}(x)=1\}})\\
&=\sum_{x\in T^d(r)}{\rm Var}(1_{\{\zeta_T^{d,p}(x)=1\}})\\
&\leq r^d P\big(\zeta_t^{d,p}(x)=1\big)=r^dpe^{-T}.
\end{align*}
As a result,
\begin{equation*}
P\big(\inf_{0\leq t\leq T}\{\frac{|A_t^{d,p}|}{r^d}-pe^{-t}\}<-\epsilon\big)
\leq \frac{e^{2T}}{\epsilon^2r^{2d}}r^dpe^{-T}=\frac{e^Tp}{\epsilon^2r^d}
\end{equation*}
and
\[
\lim_{d\rightarrow+\infty}P\big(\inf_{0\leq t\leq T}\{\frac{|A_t^{d,p}|}{r^d}-pe^{-t}\}<-\epsilon\big)=0.
\]

\qed

\section{Proof of \eqref{equ of limsup of subcri fluid}}\label{sect of proof 2.4}
In this section, we give the proof of \eqref{equ of limsup of subcri fluid}. First we introduce some notations. For any $t\geq 0$, we define
\begin{align*}
&B_t^{d,p}=\{x\in T^d(r):\eta_t^{d,p}(x)=0\},\\
&C_t^{d,p}=\{x\in T^d(r):\sum_{y:y\sim x}\eta_t^{d,p}(y)\geq d\},\\
&D_t^{d,p}=\{x\in T^d(r):\sum_{y:y\sim x}\eta_t^{d,p}(y)\leq d\},\\
&E_t^{d,p}=\bigcup_{s\leq t}C_s^{d,p}.
\end{align*}

For any $x\in A_0^{d,p}$, we define
\[
\tau_x=\inf\{t:\eta_t^{d,p}(x)=0\}
\]
and
\[
F_t^{d,p}=\{x\in A_0^{d,p}:\tau_x>t\}.
\]

The following two lemmas are crucial for our proofs.
\begin{lemma}\label{lemma 4.1}
For any $T>0$ and $t\leq T$,
\begin{equation}\label{equ (4.1)}
A_t^{d,p}\subseteq\{x\in A_0^{d,p}\setminus E_T^{d,p}:\eta_t^{d,p}(x)=1\}\cup E_T^{d,p}.
\end{equation}
\begin{equation}\label{equ (4.2)}
\{x\in A_0^{d,p}\setminus E_T^{d,p}:\eta_t^{d,p}(x)=1\}=\{x\in A_0^{d,p}\setminus E_T^{d,p}:\tau_x>t\}.
\end{equation}
\end{lemma}
\proof
\begin{equation*}
A_t^{d,p}\subseteq (A_t^{d,p}\setminus E_T^{d,p})\cup E_T^{d,p}.
\end{equation*}
For any $x\in A_t^{d,p}\setminus E_T^{d,p}$, $\eta_t^{d,p}(x)=1$ and $\sum_{y:y\sim x}\eta_t^{d,p}(y)<d$ for any $s\leq t$. According to the definition of the threshold voter model, $0$ may flip to $1$ when and only when there are not less than $d$ neighbors in state $1$. Then, if $\eta_0^{d,p}(x)=0$, $x$ can not flip to $1$ during $[0,t]$ and hence $\eta_t^{d,p}(x)=0$, which is contradictory. Therefore, $x\in A_0^{d,p}$ and
hence
\[
A_t^{d,p}\setminus E_T^{d,p}\subseteq (A_0^{d,p}\setminus E_T^{d,p})\cap A_t^{d,p},
\]
which gives \eqref{equ (4.1)}.

For any $x\in A_0^{d,p}$, if $\tau_x>t$ then $\eta_t^{d,p}(x)=1$ and hence
\[
\{x\in A_0^{d,p}\setminus E_T^{d,p}:\tau_x>t\}\subseteq \{x\in A_0^{d,p}\setminus E_T^{d,p}:\eta_t^{d,p}(x)=1\}.
\]
For any $x\in (A_0^{d,p}\setminus E_T^{d,p})\cap A_t^{d,p}$, if $\tau_x<t$, then $x$ flips from $0$ to $1$ at some moment $s\in (\tau_x,t]$. Then,
\[
\sum_{y:y\sim x}\eta_s^{d,p}(y)\geq d,
\]
which is contradictory to that $x$ does not belong to $E_T^{d,p}$.
Therefore,
\[
\{x\in A_0^{d,p}\setminus E_T^{d,p}:\eta_t^{d,p}(x)=1\}\subseteq \{x\in A_0^{d,p}\setminus E_T^{d,p}:\tau_x>t\},
\]
which gives \eqref{equ (4.2)}.

\qed

\begin{lemma}\label{lemma 4.2}
For any $T>0$, $\epsilon >0$ and $p<1/2$,
\begin{equation}\label{equ (4.3)}
\lim_{d\rightarrow+\infty}P\big(\frac{|E_T^{d,p}|}{r^d}>\epsilon\big)=0.
\end{equation}
\end{lemma}
Lemma \ref{lemma 4.2} shows that $\lim_{d\rightarrow+\infty}\frac{|E_T^{d,p}|}{r^d}=0$ in probability when $p<1/2$, which is very important for us to control $|A_t^{d,p}|$. We give the proof of Lemma \ref{lemma 4.2} at the end of this section.

Now we give the proof of \eqref{equ of limsup of subcri fluid}.

\proof[Proof of \eqref{equ of limsup of subcri fluid}]
By \eqref{equ (4.1)}, for any $t\leq T$,
\[
|A_t^{d,p}|\leq |(A_0^{d,p}\setminus E_T^{d,p})\cap A_t^{d,p}|+|E_T^{d,p}|,
\]
and hence
\[
\{\frac{|A_t^{d,p}|}{r^d}\geq pe^{-t}+\epsilon\}\subseteq \{\frac{|(A_0^{d,p}\setminus E_T^{d,p})\cap A_t^{d,p}|}{r^d}\geq pe^{-t}+\epsilon/2\} \cup \{|E_T^{d,p}|\geq \epsilon/2\}.
\]
Therefore,
\begin{align}\label{equ (4.4)}
P\big(\sup_{0\leq t\leq T}\{\frac{|A_t^{d,p}|}{r^d}-pe^{-t}\}>\epsilon\big)
\leq &P\big(\sup_{0\leq t\leq T}\{\frac{|(A_0^{d,p}\setminus E_T^{d,p})\cap A_t^{d,p}|}{r^d}-pe^{-t}\}>\epsilon/2)\notag\\
&+P\big(\frac{|E_T^{d,p}|}{r^d}>\epsilon/2\big).
\end{align}
By \eqref{equ (4.2)},
\[
|(A_0^{d,p}\setminus E_T^{d,p})\cap A_t^{d,p}|=|\{x\in A_0^{d,p}\setminus E_T^{d,p}:\tau_x>t\}|\leq |F_t^{d,p}|.
\]
Therefore,
\begin{align}\label{equ (4.5)}
&P\big(\sup_{0\leq t\leq T}\{\frac{|(A_0^{d,p}\setminus E_T^{d,p})\cap A_t^{d,p}|}{r^d}-pe^{-t}\}>\epsilon/2)\notag\\
\leq &P\big(\sup_{0\leq t\leq T}\{\frac{|F_t^{d,p}|}{r^d}-pe^{-t}\}>\epsilon/2)\leq P\big(\sup_{0\leq t\leq T}|\frac{|F_t^{d,p}|}{r^d}-pe^{-t}|>\epsilon/2\big).
\end{align}
Notice that $\{\tau_x\}_{x\in A_0^{d,p}}$ are independent and identically distributed. For any $x\in A_0^{d,p}$, $\tau_x$ is with exponential distribution with rate one. Therefore,
\[
\{F_t^{d,p}\}_{t\geq 0} \mathrel{\mathop=^{\rm d}} \{G_t^{d,p}\}_{t\geq 0}.
\]
Hence,
\begin{align}\label{equ (4.6)}
&\lim_{d\rightarrow+\infty}P\big(\sup_{0\leq t\leq T}|\frac{|F_t^{d,p}|}{r^d}-pe^{-t}|>\epsilon/2\big)\notag\\
=&\lim_{d\rightarrow+\infty}P\big(\sup_{0\leq t\leq T}|\frac{|G_t^{d,p}|}{r^d}-pe^{-t}|>\epsilon/2\big)
\leq\lim_{d\rightarrow+\infty}\frac{4e^Tp}{\epsilon^2r^d}=0
\end{align}
as we have shown in section \ref{sect of proof 2.3}.

\eqref{equ of limsup of subcri fluid} follows from \eqref{equ (4.3)}, \eqref{equ (4.4)}, \eqref{equ (4.5)}, \eqref{equ (4.6)}.

\qed

Now we only need to prove Lemma \ref{lemma 4.2}. For any integer $k\in \{0,1,2,\ldots,2d\}$ and $t>0$, we define
\begin{align*}
&H_t^{d,p}(k)=\{x\in T^d(r):\sum_{y:y\sim x}\eta_t^{d,p}(y)=k\},\\
&I_t^{d,p}(k)=\{x\in T^d(r):\sum_{y:y\sim x}\eta_t^{d,p}(y)\geq k\}=\bigcup_{l\geq k}H_t^{d,p}(l),\\
&J_t^{d,p}(k)=\{x\in T^d(r):\sum_{y:y\sim x}\eta_t^{d,p}(y)\leq k\}=\bigcup_{l\leq k}H_t^{d,p}(l).
\end{align*}
We introduce two lemmas to prove Lemma \ref{lemma 4.2}.
\begin{lemma}\label{lemma 4.3}
Assume that $p<1/2$, then for any $q\in (p,1]$, $\epsilon>0$,
\[
\lim_{d\rightarrow+\infty}P\big(\frac{|I_0^{d,p}(\lfloor2dq\rfloor)|}{r^d}>\epsilon\big)=0.
\]
\end{lemma}
\proof
\begin{align}\label{equ (4.7)}
{\rm E}|I_0^{d,p}(\lfloor2dq\rfloor)|&={\rm E}\sum_x1_{\{\sum_{y:y\sim x}\eta_0^{d,p}(y)\geq \lfloor2dq\rfloor\}}\notag\\
&=\sum_xP\big(\sum_{y:y\sim x}\eta_0^{d,p}(y)\geq\lfloor2dq\rfloor\big)\notag\\
&=r^dP\big(\sum_{y:y\sim O}\eta_0^{d,p}(y)\geq\lfloor2dq\rfloor\big),
\end{align}
where $O$ is a fixed vertex on $T^d(r)$.

Since $\eta_0^{d,p}$ is with distribution $\mu_p$, $\{\eta_0^{d,p}(y)\}_{y:y\sim O}$ are independent and identically distributed with
\[
P\big(\eta_0^{d,p}(y)=1\big)=p=1-P\big(\eta_0^{d,p}(y)=0\big).
\]
Therefore, according to law of large numbers,
\begin{equation}\label{equ (4.8)}
\lim_{d\rightarrow+\infty} \frac{\sum_{y:y\sim O}\eta_0^{d,p}(y)}{2d}=p
\end{equation}
in probability.

By \eqref{equ (4.7)}, \eqref{equ (4.8)} and Chebyshev's inequality,
\begin{align*}
&\lim_{d\rightarrow+\infty}P\big(\frac{|I_0^{d,p}(\lfloor2dq\rfloor)|}{r^d}>\epsilon\big)\\
\leq&\lim_{d\rightarrow+\infty}\frac{1}{r^d\epsilon}{\rm E}|I_0^{d,p}(\lfloor2dq\rfloor)|\\
=&\lim_{d\rightarrow+\infty}\frac{1}{\epsilon}P\big(\sum_{y:y\sim O}\eta_0^{d,p}(y)\geq\lfloor2dq\rfloor\big)\\
\leq&\lim_{d\rightarrow+\infty}\frac{1}{\epsilon}P\big(\frac{\sum_{y:y\sim O}\eta_0^{d,p}(y)}{2d}-p>0.9(q-p)\big)=0.
\end{align*}

\qed

\begin{lemma}\label{lemma 4.4}
There exists $\epsilon_r>0$ such that for any $p\in (1/2-\epsilon_r,1/2)$, there exists $C(p)>0$ such that
\[
\lim_{d\rightarrow+\infty}\frac{|C_0^{d,p}|}{e^{C(p)d}}=+\infty
\]
in probability.
\end{lemma}
\proof
\[
{\rm E}|C_0^{d,p}|=r^dP\big(\sum_{y:y\sim O}\eta^{d,p}_0(y)\geq d\big).
\]
According to Cram\'{e}r's Theorem (See Chapter two of \cite{Dem1997}),
\[
\lim_{d\rightarrow+\infty}\frac{1}{d}\log P\big(\sum_{y:y\sim O}\eta^{d,p}_0(y)\geq d\big)=-K(p),
\]
where
\[
K(p)=-\log[4p(1-p)]>0
\]
when $p<1/2$. $K(1/2)=0$ and $K(p)$ is continuous with $p$, so we can choose sufficiently small $\epsilon_r$ such that
\[
\log r-K(p)>0
\]
when $p\in (1/2-\epsilon_r,1/2)$.

For $p\in (1/2-\epsilon_r,1/2)$, we let
\[
C(p)=\frac{\log r-K(p)}{2}>0.
\]
Now we only need to show that for any $M>0$,
\begin{equation}\label{equ (4.9)}
\lim_{d\rightarrow+\infty}P\big(\frac{|C_0^{d,p}|}{e^{C(p)d}}<M\big)=0.
\end{equation}
When $p\in (1/2-\epsilon_r,1/2)$,
\[
{\rm E}|C_0^{d,p}|=r^dP\big(\sum_{y:y\sim O}\eta^{d,p}_0(y)\geq d\big)=\exp\big\{2dC(p)+o(d)\big\}.
\]

Then according to Chebyshev's inequality, for sufficiently large $d$,
\begin{equation}\label{equ (4.10)}
P\big(\frac{|C_0^{d,p}|}{e^{C(p)d}}<M\big)<\frac{{\rm Var}(|C_0^{d,p}|)}{\exp\big\{4dC(p)+o(d)\big\}}.
\end{equation}
Now we calculate ${\rm Var}(|C_0^{d,p}|)$,
\begin{align}\label{equ (4.11)}
&{\rm Var}(|C_0^{d,p}|)
={\rm E}|C_0^{d,p}|^2-(E|C_0^{d,p}|)^2\notag\\
=&{\rm E}(\sum_{x\in T^d(r)}1_{\{\sum_{y:y\sim x}\eta_0^{d,p}(y)\geq d\}})^2-[r^dP\big(\sum_{y:y\sim O}\eta^{d,p}_0(y)\geq d\big)]^2\notag\\
=&\sum_{x_1\in T^d(r)}\sum_{x_2\in T^d(r)}P\big(\sum_{y:y\sim x_1}\eta_0^{d,p}(y)\geq d,\sum_{z:z\sim x_2}\eta_0^{d,p}(z)\geq d\big)\\
&-[r^dP\big(\sum_{y:y\sim O}\eta^{d,p}_0(y)\geq d\big)]^2.\notag
\end{align}
For any $x_1,x_2\in T^d(r)$, we define
\[
N(x_1,x_2)=\{y\in T^d(r):y\sim x_1, y\sim x_2\}.
\]
According to the structure of a torus, $|N(x,y)|\leq 2$ and
\begin{equation}\label{equ (4.12)}
|\{z\in T^d(r):N(z,x)\neq\emptyset\}|=2d+4{d \choose 2}=2d^2
\end{equation}
for any $x,y\in T^d(r)$.

If $N(x_1,x_2)=\emptyset$, then
\begin{equation}\label{equ (4.13)}
P\big(\sum_{z:z\sim x_1}\eta_0^{d,p}(z)\geq d,\sum_{w:w\sim x_2}\eta_0^{d,p}(w)\geq d\big)=[P\big(\sum_{y:y\sim O}\eta^{d,p}_0(y)\geq d\big)]^2.
\end{equation}
If $N(x_1,x_2)\neq\emptyset$, then by $|N(x_1,x_2)|\leq 2$ and Cram\'{e}r's Theorem,
\begin{align}\label{equ (4.14)}
&P\big(\sum_{z:z\sim x_1}\eta_0^{d,p}(z)\geq d,\sum_{w:w\sim x_2}\eta_0^{d,p}(w)\geq d\big)\notag\\
\leq&P\big(\sum_{z:z\sim x_1,\atop z\not\in N(x_1,x_2)}\eta_0^{d,p}(z)\geq d-2,\sum_{w:w\sim x_2,\atop w\not\in N(x_1,x_2)}\eta_0^{d,p}(w)\geq d-2\big)\notag\\
=&[P\big(\sum_{z:z\sim x_1,\atop z\not\in N(x_1,x_2)}\eta_0^{d,p}(z)\geq d-2\big)]^2\notag\\
\leq&[P\big(\sum_{y:y\sim O}\eta^{d,p}_0(y)\geq d-2\big)]^2\notag\\
=&\exp\{-2dK(p)+o(d)\}.
\end{align}
By \eqref{equ (4.12)}, \eqref{equ (4.13)} and \eqref{equ (4.14)},
\begin{align}\label{equ (4.15)}
&\sum_{x_1\in T^d(r)}\sum_{x_2\in T^d(r)}P\big(\sum_{y:y\sim x_1}\eta_0^{d,p}(y)\geq d,\sum_{z:z\sim x_2}\eta_0^{d,p}(z)\geq d\big)\notag\\
\leq&r^{2d}[P\big(\sum_{y:y\sim O}\eta^{d,p}_0(y)\geq d\big)]^2+2r^dd^2\exp\{-2dK(p)+o(d)\}\notag\\
=&[r^dP\big(\sum_{y:y\sim O}\eta^{d,p}_0(y)\geq d\big)]^2+\exp\{d(\log r-2K(p))+o(d)\}.
\end{align}
By \eqref{equ (4.11)} and \eqref{equ (4.15)},
\begin{equation}\label{equ (4.16)}
{\rm Var}(|C_0^{d,p}|)\leq \exp\{d(\log r-2K(p))+o(d)\}.
\end{equation}
Then by \eqref{equ (4.10)},
\[
P\big(\frac{|C_0^{d,p}|}{e^{C(p)d}}<M\big)<\exp\{-d\log r+o(d)\}
\]
and hence
\[
\lim_{d\rightarrow+\infty}P\big(\frac{|C_0^{d,p}|}{e^{C(p)d}}<M\big)=0.
\]

\qed

Now we give the proof of Lemma \ref{lemma 4.2}. The proof is inspired a lot by the approach of moving balls between boxes introduced in \cite{Lanchier2010}.
\proof[Proof of Lemma \ref{lemma 4.2}]
It is easy to see that the threshold voter model is an attractive spin system (See Section 3.1 of \cite{LIG1985}). Hence when $p_1<p_2<1/2$,
\[
P\big(\frac{|E_T^{d,p_1}|}{r^d}>\epsilon\big)\leq P\big(\frac{|E_T^{d,p_2}|}{r^d}>\epsilon\big).
\]
So we only need to deal with the case that $p\in (1/2-\epsilon_r,1/2)$.

According to the flip rates given by \eqref{equ of flip rate of threshold voter on tori}, $|A_t^{d,p}|$ evolves as follows,
\begin{equation}
|A_t^{d,p}|\rightarrow
\begin{cases}
|A_t^{d,p}|+1 \text{~at rate $|B_t^{d,p}\cap C_t^{d,p}|$,}\\
|A_t^{d,p}|-1 \text{~at rate $|A_t^{d,p}\cap D_t^{d,p}|$.}
\end{cases}
\end{equation}
We set $2d+1$ boxes $b_0,b_1,\ldots, b_{2d}$. At $t=0$, we put $|H_0^{d,p}(k)|$ balls in box $b_k$ for $k\in \{0,1,2,\ldots,2d\}$. For each $k$, we say $b_{k-1}$ is the left box of $b_k$ and $b_{k+1}$ is the right box of $b_k$.

Now we introduce four approaches of moving balls between boxes. In the first approach,  we move balls according to the evolution of $\eta_t^{d,p}$ such that at any $t>0$, there are $|H_t^{d,p}(k)|$ balls in box $b_k$ for each $k$. Each vertex in $T^d(r)$ has $2d$ neighbors, therefore, at rate $|B_t^{d,p}\cap C_t^{d,p}|$, $2d$ balls move to their right boxes. At rate $|A_t\cap D_t^{d,p}|$, $2d$ balls move to their left boxes. $|C_t^{d,p}|$ is the sum of the numbers of balls in boxes $\{b_k\}_{d\leq k\leq 2d}$.

To control $|C_t^{d,p}|$ from above, we introduce a second approach of moving balls. We denote by $b_k(t)$ the number of balls in $b_k$ at $t$ and define
\[
\widehat{C}_t=\sum_{k=d}^{2d}b_k(t).
\]
In this approach, we never move balls to their left boxes. For any $t\geq0$, we move $2d$ balls to right boxes at rate $\widehat{C}_t$. These $2d$ balls are chosen and moved as follows.
Let
\[
j_1=\sup\{k<d:b_k(t)>0\}
\]
and
\[
j_l=\sup\{k<j_{l-1}:b_k(t)>0\}
\]
for $l=2,3,\ldots.$

Let
\[
\tau=\inf\{k:\sum_{l=1}^kb_{j_l}(t)\geq 2d\}.
\]
If $\tau<+\infty$, then we move $b_{j_l}(t)$ balls from $b_{j_l}$ to $b_{j_l+1}$ for each $l<\tau$ and move $2d-\sum_{l=1}^{\tau-1}b_{j_l}(t)$ balls from $b_{j_\tau}$
to $b_{j_\tau+1}$. If $\tau=+\infty$, then we move $b_l(t)$ balls from $b_l$ to $b_{l+1}$ for each $l<d$ and pick $2d-\sum_{l=0}^{d-1}b_l(t)$ balls randomly from boxes $\{b_k\}_{d\leq k\leq 2d}$ and put these $2d-\sum_{l=0}^{d-1}b_l(t)$ balls in $b_{2d}$.

In the second approach, we choose $2d$ balls which are nearest to $b_d$ to move right. This approach makes the total number of balls in boxes $\{b_k\}_{d\leq k\leq 2d}$ increase as fast as possible. We never move balls to left boxes, so when a ball reaches $b_k$ with $k\geq d$, it will never return to any $b_j$ with $j<d$. As a result, for any $M>0$,
\begin{equation}\label{equ (4.18)}
P\big(|E_T^{d,p}|>M\big)=P\big(|\cup_{s\leq T}C_s^{d,p}|>M\big)\leq P\big(\widehat{C}_T>M\big).
\end{equation}

We modify the second approach at $t=0$ to obtain the third approach. Let $p_0=1/4+p/2\in (p,1/2)$. At $t=0$, we pick all the balls in boxes $\{b_k:\lfloor2dp_0\rfloor\leq k\leq d-1\}$ and put these balls in $b_d$. Then, we still move balls as the way in the second approach. We denote by $\overline{C}_t$ the total number of balls in boxes $\{b_k\}_{d\leq k\leq 2d}$ at $t$. The third approach accelerate the moving-ball process in the second approach, since the original several steps of moving balls in boxes $\{b_k:\lfloor2dp_0\rfloor\leq k\leq d-1\}$ are finished at $0$.

Therefore,
\begin{equation}\label{equ (4.19)}
P\big(\widehat{C}_T>M\big)\leq P\big(\overline{C}_T>M\big).
\end{equation}

Notice that in the third approach, $\overline{C}_0=|I_0^{d,p}(\lfloor2dp_0\rfloor)|$ and there are no balls in $b_k$ for $\lfloor 2dp_0\rfloor\leq k<d$ at $t=0$. Then it takes not less than $\lfloor d(1-2p_0)\rfloor$ steps of moving balls to make $\overline{C}$ add $2d$. Inspired by this phenomenon, we construct the fourth approach. In the fourth approach, there is only one box. At $t=0$, we put $|I_0^{d,p}(\lfloor2dp_0\rfloor)|$ balls in the box. $\{Y_i\}_{i\geq 1}$ are independent and identically distributed random variables with exponential distribution with rate $1$. Further more, we assume that $\{Y_i\}_{i\geq 1}$ are independent with $\eta_0^{d,p}$. For each $j\geq 1$, we define
\[
\tau_j=\sum_{l=(j-1)\lfloor d(1-2p_0)\rfloor+1}^{j\lfloor d(1-2p_0)\rfloor}\frac{Y_l}{|I_0^{d,p}(\lfloor2dp_0\rfloor)|+2d(j-1)}
\]
and
\[
t_j=\sum_{l=1}^j\tau_l.
\]
We add $2d$ balls in the box at $t_j$ for each $j\geq 1$. For $t\not \in \{t_j\}_{j\geq 1}$, the balls' number stay still. We denote by $\widetilde{C}_t^{d,p}$ the balls' number at $t$. Then
\begin{equation*}
\widetilde{C}_t^{d,p}=|I_0^{d,p}(\lfloor2dp_0\rfloor)|+2dj
\end{equation*}
for $t\in [t_j,t_{j+1})$.

Notice that $\frac{Y_l}{\alpha}$ is with exponential distribution with rate $\alpha$ for any $\alpha>0$. So $\tau_j$ is the sum of $\lfloor d(1-2p_0)\rfloor$ i.i.d random variables with exponential distribution with rate $\widetilde{C}_0^{d,p}+2d(j-1)$. The number of balls in the box
takes $\tau_j$ to increase to $\widetilde{C}_0^{d,p}+2dj$ from $\widetilde{C}_0^{d,p}+2d(j-1)$. Compared the third and fourth approaches,
\begin{equation}\label{equ (4.20)}
P\big(\overline{C}_T>M\big)\leq P\big(\widetilde{C}_T^{d,p}>M\big).
\end{equation}
By \eqref{equ (4.18)}, \eqref{equ (4.19)} and \eqref{equ (4.20)}, we only need to show that
\begin{equation}\label{equ (4.21)}
\lim_{d\rightarrow+\infty}P\big(\frac{\widetilde{C}_T^{d,p}}{r^d}>\epsilon\big)=0
\end{equation}
for any $\epsilon$.

Notice that,
\begin{align*}
P\big(\frac{\widetilde{C}_T^{d,p}}{r^d}\leq\epsilon\big)
&\geq P\big(\widetilde{C}_T^{d,p}\leq\widetilde{C}_0^{d,p}e^{\frac{4T}{1-2p_0}},\frac{\widetilde{C}_0^{d,p}}{r^d}\leq
\epsilon e^{-\frac{4T}{1-2p_0}}\big)\\
&=P\big(\widetilde{C}_T^{d,p}\leq\widetilde{C}_0^{d,p}e^{\frac{4T}{1-2p_0}},\frac{|I_0^{d,p}(\lfloor2dp_0\rfloor)|}{r^d}\leq
\epsilon e^{-\frac{4T}{1-2p_0}}\big).
\end{align*}
By Lemma \ref{lemma 4.3},
\[
\lim_{d\rightarrow+\infty}P\big(\frac{|I_0^{d,p}(\lfloor2dp_0\rfloor)|}{r^d}\leq
\epsilon e^{-\frac{4T}{1-2p_0}}\big)=1.
\]
Therefore, to prove \eqref{equ (4.21)}, we only need to show that
\begin{equation}\label{equ (4.22)}
\lim_{d\rightarrow+\infty}P\big(\widetilde{C}_T^{d,p}\geq\widetilde{C}_0^{d,p}e^{\frac{4T}{1-2p_0}}\big)=0.
\end{equation}
According to Lemma \ref{lemma 4.4}, for $p\in (1/2-\epsilon_r,1/2)$,
\begin{equation}\label{equ (4.23)}
\lim_{d\rightarrow+\infty}P\big(\widetilde{C}_0^{d,p}\leq e^{C(p)d}\big)=0
\end{equation}
for some $C(p)>0$.
\begin{align}\label{equ (4.24)}
P\big(\widetilde{C}_T^{d,p}\geq\widetilde{C}_0^{d,p}e^{\frac{4T}{1-2p_0}}\big)
\leq &P\big(\widetilde{C}_T^{d,p}\geq\widetilde{C}_0^{d,p}e^{\frac{4T}{1-2p_0}}\big|\widetilde{C}_0^{d,p}>e^{C(p)d}\big)\notag\\
&+P\big(\widetilde{C}_0^{d,p}\leq e^{C(p)d}\big).
\end{align}
To distinguish the dimension $d$, we write $t_j$ and $\tau_j$ as $t_j^d$ and $\tau_j^d$.
Let
\[
K_d=\lfloor\frac{\widetilde{C}_0^{d,p}(e^{\frac{4T}{1-2p_0}}-1)}{2d}\rfloor.
\]
Conditioned on $\widetilde{C}_0^{d,p}>e^{C(p)d}$ for each $d$,
\begin{equation}
\lim_{d\rightarrow+\infty}K_d=+\infty
\end{equation}
and
\begin{equation}\label{equ (4.26)}
P\big(\widetilde{C}_T^{d,p}\geq\widetilde{C}_0^{d,p}e^{\frac{4T}{1-2p_0}}\big)
\leq P\big(t^d_{K_d}\leq T\big)=P\big(\sum_{l=1}^{K_d}\tau^d_l\leq T\big).
\end{equation}
For $j\geq 1$, conditioned on $\widetilde{C}_0^{d,p}$,
\[
{\rm E}\tau_j^d=\frac{\lfloor d(1-2p_0)\rfloor}{\widetilde{C}_0^{d,p}+2d(j-1)}.
\]
We choose sufficiently small $\delta$ such that $\delta<0.1$ and
\[
\log\frac{e^{\frac{4T}{1-2p_0}}}{1-\delta+\delta e^{\frac{4T}{1-2p_0}}}>\frac{3.9T}{1-2p_0}.
\]

For $j\geq 1$, let
\[
Z_j^d=
\begin{cases}
1\text{~if $\frac{\tau_j^d}{{\rm E}\tau_j^d}>1-\delta$,}\\
0\text{~else.}
\end{cases}
\]
If
\[
\frac{\sum_{l=1}^{K_d}Z_l^d}{K_d}>1-\delta,
\]
then
\begin{align*}
\sum_{l=1}^{K_d}\tau_l^d&\geq\sum_{l=\lceil\delta K_d\rceil}^{K_d}\frac{0.9d(1-2p_0)}{\widetilde{C}_0^{d,p}+2d(l-1)}+o(1)\\
&=(0.45-0.9p_0)\sum_{l=\lceil\delta K_d\rceil}^{K_d}\frac{2d}{\widetilde{C}_0^{d,p}+2d(l-1)}+o(1)\\
&\geq (0.45-0.9p_0)\int_{\widetilde{C}_0^{d,p}+2d\delta K_d}^{\widetilde{C}_0^{d,p}+2dK_d}\frac{1}{t}dt+o(1)\\
&=(0.45-0.9p_0)\log\frac{e^{\frac{4T}{1-2p_0}}}{1-\delta+\delta e^{\frac{4T}{1-2p_0}}}+o(1)\\
&\geq 1.75T
\end{align*}
for sufficiently large $d$.

Therefore, for sufficiently large $d$,
\begin{equation}\label{equ (4.27)}
P\big(\sum_{l=1}^{K_d}\tau^d_l\leq T\big)\leq P\big(\frac{\sum_{l=1}^{K_d}Z_l^d}{K_d}\leq1-\delta\big).
\end{equation}
By the definition of $\tau_j$,
\[
\frac{\tau_j^d}{{\rm E}\tau_j^d}\mathrel{\mathop=^{\rm d}}\frac{1}{\lfloor d(1-2p_0)\rfloor}\sum_{l=1}^{\lfloor d(1-2p_0)\rfloor}Y_l
\]
for each $j\geq 1$.

Therefore, according to the law of large numbers,
\begin{equation}\label{equ (4.28)}
P\big(\frac{\tau_j^d}{{\rm E}\tau_j^d}>1-\delta\big)>1-\frac{\delta}{10}
\end{equation}
for sufficiently large $d$.

Let $\{W_i\}_{i\geq 1}$ be independent and identically distributed random variables with
\[
P\big(W_1=1\big)=1-\frac{\delta}{10}=1-P\big(W_1=0\big).
\]
By \eqref{equ (4.28)}, for sufficiently large $d$,
\begin{equation}\label{equ (4.29)}
P\big(\frac{\sum_{l=1}^{K_d}Z_l^d}{K_d}\leq1-\delta\big)\leq P\big(\frac{\sum_{l=1}^{K_d}W_l}{K_d}\leq 1-\delta\big).
\end{equation}
Since $K_d\rightarrow+\infty$, according to law of large numbers,
\[
\lim_{d\rightarrow+\infty}P\big(\frac{\sum_{l=1}^{K_d}W_l}{K_d}\leq 1-\delta\big)=0
\]
and
\begin{equation}\label{equ (4.30)}
\lim_{d\rightarrow+\infty}P\big(\frac{\sum_{l=1}^{K_d}Z_l^d}{K_d}\leq1-\delta\big)=0.
\end{equation}
\eqref{equ (4.22)} follows from \eqref{equ (4.23)}, \eqref{equ (4.24)}, \eqref{equ (4.26)}, \eqref{equ (4.27)} and \eqref{equ (4.30)}.
We have shown that Lemma \ref{lemma 4.2} follows from \eqref{equ (4.22)}. The proof is finished.

\qed

\textbf{Acknowledgments.} This work is supported by the National Basic Research Program of China (2011CB808000).

{}
\end{document}